\documentclass[12pt,a4paper]{amsart}
\usepackage{amsmath}
\usepackage{cases}
\usepackage{ifthen,verbatim}
\usepackage{graphicx}
\usepackage{mathrsfs}
\nonstopmode
\numberwithin{equation}{section}
\newtheorem{thm}{Theorem}[section]

\newtheorem{lem}[thm]{Lemma}
\newtheorem{prop}[thm]{Proposition}

\theoremstyle{definition}

\newtheorem*{ThmA}{Theorem A}
\newtheorem*{ThmB}{Theorem B}
\newtheorem*{ThmC}{Theorem C}

\newenvironment{pf}[1][]{%
 \vskip 3mm
 \noindent
 \ifthenelse{\equal{#1}{}}%
  {{\slshape Proof. }}%
  {{\slshape #1.} }%
 }%
{\qed\bigskip}

\newcounter{alphabet}
\newcounter{tmp}

\usepackage{cases}

\newcommand{\B}{{\mathscr B}}
\newcommand{\C}{{\mathbb C}}
\newcommand{\D}{{\mathbb D}}
\newcommand{\F}{{\mathscr F}}

\newcommand{\PP}{{\mathscr P}}
\newcommand{\N}{{\mathbb N}}
\newcommand{\R}{{\mathbb R}}

\newcommand{\Schur}{{\mathscr{S}}}
\newcommand{\T}{{\mathbb T}}
\newcommand{\U}{{\mathbf U}}
\newcommand{\V}{{\mathbf V}}
\newcommand{\X}{{\mathbf X}}

\newcommand{\hol}{{\operatorname{Hol}}}
\newcommand{\Int}{{\operatorname{Int}\,}}

\newcommand{\bD}{{\overline{\mathbb D}}}

\renewcommand{\Re}{{\,\operatorname{Re}\,}}
\newcommand{\id}{{\operatorname{id}}}

\newcommand{\aand}{{\quad\text{and}\quad}}

\begin{document}
\bibliographystyle{amsplain}

\title[Schur parameters and Carath\'eodory class]%
{Schur parameters and Carath\'eodory class}
\begin{abstract}
The Schur (resp.~ Carath\'eodory) class consists of all the analytic functions $f$ on the unit disk with
$|f|\le 1$ (resp.~$\Re f>0$ and $f(0)=1$).
The Schur parameters $\gamma_0,\gamma_1,\dots (|\gamma_j|\le 1)$ are known to
parametrize the coefficients of functions in the Schur class.
By employing a recursive formula for it,
we describe the $n$-th coefficient of a Carath\'eodory function
in terms of $n$ independent variables $\gamma_1,\dots, \gamma_n$ with $|\gamma_j|\le 1.$
The mapping properties of those correspondences are also studied.
\end{abstract}
\date{\today}

\author[M. Li]{Ming Li}
\address{Institute for Advanced Study \\
Shenzhen University \\
Nanhai Ave 3688, Shenzhen, Guangdong 518060, P.~R.~China}
\email{mingli@szu.edu.cn / minglimath@sina.com}
\author[T. Sugawa]{Toshiyuki Sugawa}
\address{Graduate School of Information Sciences \\
Tohoku University \\
Aoba-ku, Sendai 980-8579, Japan}
\email{sugawa@math.is.tohoku.ac.jp}

\maketitle

\section{Introduction }
We denote by $\hol(\D)$ the set of analytic ($=$ holomorphic) functions on
the unit disk $\D=\{z: |z|<1\}$ in the complex plane $\C.$
In the present note, our main concern is about the subclasses
\begin{align*}
\Schur&=\{f\in\hol(\D): |f|\le 1\}, \\
\PP&=\{g\in\hol(\D): g(0)=1,~ \Re g>0\}.
\end{align*}
These are called the Schur class and the Carath\'eodory class, respectively, and
members of $\Schur$ and $\PP$ are called Schur functions and Carath\'eodory
functions, respectively.
We also consider the subclass $\Schur_0=\{\omega\in\Schur: \omega(0)=0\}
=\{zf(z): f\in\Schur\}$ of $\Schur.$ 
These classes play an important role in Geometric Function Theory.
For example, a function $f\in\hol(\D)$ with $f(0)=0, f'(0)=1,$
is starlike (resp.~ convex) if $zf'(z)/f(z)$ (resp.~$1+zf''(z)/f'(z)$) belongs
to the class $\PP.$
Also, a function $f$ in $\hol(\D)$ is said to be subordinate to another $g\in\hol(\D)$
(we write $f\prec g$ for it) if $f=g\circ\omega$ for some $\omega\in\Schur_0.$
Therefore, detailed information about the Taylor coefficients of functions in
$\Schur$ and $\PP$ will lead to various useful estimates in Geometric Function Theory.
As is well known, a function $f(z)=c_0+c_1z+c_2z^2+\cdots$ in $\Schur$ satisfies $|c_n|\le 1$
for all $n\ge 0$ and the bound is sharp for each $n.$
For a function $g(z)=1+p_1z+p_2z^2+\cdots$ in $\PP,$ the sharp inequality
$|p_n|\le 2~(n=1,2,\dots)$ is known as the Carath\'eodory lemma.
A complete characterization of the coefficients of $\Schur$ and $\PP$ in terms of determinants
are classically known.
For instance, the coefficients of Carath\'eodory functions are described in the following theorem
due to Carath\'eodory and Toeplitz (see \cite{Tsuji:Potential}).

\begin{ThmA}\label{Thm:CT}
Let $g(z)=1+p_1z+p_2z^2+\cdots$ be a formal power series with coefficients in $\C.$
Then, $g$ represents a Carath\'eodory function on $\D$ if and only if
\begin{equation*}
\Delta_n=
\begin{vmatrix}
2 & p_1& p_2 &\cdots& p_n\\
\overline{p_1} & 2  & p_1 &\cdots & p_{n-1}\\
\overline{p_2} & \overline{p_1} & 2  &\cdots & p_{n-2}\\
\vdots&\vdots&\vdots&\ddots&\vdots\\
\overline{p_n} & \overline{p_{n-1}} & \overline{p_{n-2}} &\cdots & 2
\end{vmatrix}\ge0
\end{equation*}
for all $n\ge0.$
\end{ThmA}

Therefore, the region of the coefficients $p_1,\dots, p_n$ is described by the $n$ inequalities
$\Delta_1\ge0,\dots, \Delta_n\ge0.$
For instance, the region of the first coefficient is $|p_1|\le 2.$
However, it is not easy to use them in general.
Parametric representations of the coefficients are often more useful.
Libera and Z\l otkiewicz \cite{LZ82} derived the following parametrizations
of possible values of $p_2 $ and $p_3$ from Theorem A.

\begin{ThmB}\label{Thm:LZ}
Let $g(z)=1+p_1z+p_2z^2+\cdots$ be a Carath\'eodory function with $p_1\ge0.$
Then there are numbers $x, y\in\bD$ such that
\begin{align*}
2p_2&=p_1^2+x(4-p_1^2) \quad\text{and} \\
4p_3&=p_1^3+(4-p_1^2)(2p_1x -p_1x^2)+2(4-p_1^2)(1-|x|^2)y.
\end{align*}
\end{ThmB}

Also, recently Kwon, Lecko and Sim \cite{KLS18} obtained a similar parametrization of $p_4.$
We remark that the assumption $p_1\ge0$ is harmless because
we can normalize any function $g\in\PP$ so that $p_1\ge0$ 
by considering a suitable rotation $g(e^{i\theta}z).$
In recent years, these results are used frequently to estimate Hankel determinants
of functions in special classes of univalent functions. 
See, for instance,  \cite{KLS18B} and \cite{Zap18} and references therein.
On the other hand, it is known that the coefficients of a function in the Schur class $\Schur$
are described by the Schur parameters (see \cite[\S 1.3]{Simon:OP1}).
In Section 2, we will give a recursive method to compute $c_n$ in terms of the Schur parameters,
as well as basic facts about the Schur class.
Since the classes $\Schur_0$ and $\PP$ are related by the Cayley transformation,
we can develop a systematic approach to get concrete relations between the coefficients of these
two classes in Section 3.
Then in Section 4, we will parametrize $b_n=p_n/2$ as a function $T_n(\gamma_1,\dots,\gamma_n)$
of $n$ independent variables $\gamma_1, \dots, \gamma_n$ in $\bD.$
The definition of $T_n$ given in Section 4 is convenient to observe basic properties of it
but somewhat indirect.
We also give a recursive formula to describe $T_n$ in Section 5.

\section{Schur algorithm}

Let $a\in\D.$
We first recall that the M\"obius transformation
$$
L_a(z)=\frac{z-a}{1-\bar a z}
$$
keeps $\D$ invariant as a set and maps $a$ to $0$ and $0$ to $-a.$
In particular, $L_a\in\Schur.$
A function in $\Schur$ of the form
$$
e^{i\theta}L_{a_1}(z)\cdots L_{a_d}(z)
\quad (\theta\in\R, a_1,\dots, a_d\in\D)
$$
is called a finite Blaschke product of degree $d$ and the set of all
such functions will be denoted by $\B_d.$
Note that $\B_0$ consists of unimodular constants; namely,
$\B_0=\T=\partial\D,$
and that $\B_1$ is nothing but the group of analytic automorphisms of $\D.$
For a function $f\in\Schur\setminus\B_0,$ we consider the new function $\sigma f$
defined by
$$
\sigma f(z)=\frac1z\cdot \frac{f(z)-f(0)}{1-\overline{f(0)}f(z)}=\frac{L_{f(0)}(f(z))}{z}.
$$
Since the origin is a removable singularity, the function $\sigma f$ is analytic on $\D.$
Moreover, the maximum modulus principle implies that 
$$
|\sigma f(z)|\le \lim_{r\to1^-}\max_{|\zeta|\le r}|\sigma f(\zeta)|\le \lim_{r\to1^-}\frac1r=1,
$$
and hence $\sigma f\in\Schur.$
We define $\sigma f=0$ when $f\in\B_0.$
In this way, a self-map $\sigma:\Schur\to\Schur$ is defined.
For a given $f\in\Schur,$ we start with $f_0=f$ and define $f_n$
inductively by $f_n=\sigma f_{n-1}$ for $n\ge 1.$
That is to say, $f_n=\sigma^n f.$
This procedure is usually called the Schur algorithm.
We define a sequence $\gamma_n~(n=0,1,2,\dots)$ by setting $\gamma_n=f_n(0)$
and call those numbers {\it Schur parameters}.
For convenience, sometimes we write $\vec{\gamma}=\vec{\gamma}(f)=(\gamma_0,\gamma_1,\dots)
\in \bD^{\N_0},$ where $\N_0=\{0,1,2,\dots\},$ and call it the Schur vector of $f.$
By definition, we observe that $\vec{\gamma}(\sigma f)=(\gamma_1,\gamma_2,\dots)$ if
$\vec{\gamma}(f)=(\gamma_0,\gamma_1,\dots).$
What is the same, the Schur vector of $\sigma f$ is the backward shift of the Schur vector of $f.$
Schur \cite{Schur17} 
proved that the original function $f$ can be reproduced by its Schur vector.

\begin{ThmC}\label{Thm:Schur}
For a function $f\in\Schur,$ its Schur parameters $\gamma_0,\gamma_1,\dots$ satisfy
one of the following two conditions:
\begin{enumerate}
\item[(i)]
$|\gamma_n|<1$ for all $n.$
\item[(ii)]
$|\gamma_0|<1,\dots, |\gamma_{n-1}|<1, |\gamma_n|=1,
\gamma_{n+1}=0, \gamma_{n+2}=0,\dots$ for some $n\ge 0.$
\end{enumerate}
The latter occurs if and only if $f\in\B_n.$
Moreover, for any sequence $\vec{\gamma}=(\gamma_0,\gamma_1,\dots)$
satisfying one of the above conditions, there exists a unique function
$f\in\Schur$ such that $\vec{\gamma}(f)=\vec{\gamma}.$
\end{ThmC}

We remark that the subclass $\Schur_0$ is characterized by $\gamma_0=0$
with the Schur parameters.

We define a sequence of functions $F_n$ of the $n$ complex variables
$\gamma_1, \gamma_2, \dots, \gamma_n~(n=1,2,3,\dots) $ inductively by
\begin{align}
\notag
F_1(\gamma_1)&=\gamma_1, \\
\label{eq:F}
F_{n}(\gamma_1,\dots, \gamma_{n})
&=(1-|\gamma_1|^2)F_{n-1}(\gamma_2,\dots, \gamma_{n}) \\
\notag
&\quad -\bar\gamma_1\sum_{k=2}^{n-1} F_{n-k}(\gamma_2,\dots,\gamma_{n-k+1})
F_k(\gamma_1,\dots, \gamma_k), \qquad n\ge 2.
\end{align}
By construction, we see that the function $F_n(\gamma_1,\dots,\gamma_n)$
is indeed a polynomial in $\gamma_1,\bar\gamma_1, \dots,$
$\gamma_{n-1},\bar\gamma_{n-1}, \gamma_n$
with integer coefficients.
We compute the first several as follows:
\begin{align*}
F_1&=\gamma_1, \\
F_2&=(1-|\gamma_1|^2)\gamma_2, \\
F_3&=(1-|\gamma_1|^2)(1-|\gamma_2|^2)\gamma_3
-(1-|\gamma_1|^2)\bar\gamma_1\gamma_2^2, \\
F_4&=(1-|\gamma_1|^2)(1-|\gamma_2|^2)(1-|\gamma_3|^2)\gamma_4 \\
&\quad -(1-|\gamma_1|^2)(1-|\gamma_2|^2)\gamma_3(2\bar\gamma_1\gamma_2
+\bar\gamma_2\gamma_3)+(1-|\gamma_1|^2)\bar\gamma_1^2\gamma_2^3, \\
F_5&=(1-|\gamma_1|^2)(1-|\gamma_2|^2)(1-|\gamma_3|^2)(1-|\gamma_4|^2)\gamma_5 \\
&\quad  -(1-|\gamma_1|^2)(1-|\gamma_2|^2)(1-|\gamma_3|^2)\gamma_4
(2\bar\gamma_1\gamma_2+2\bar\gamma_2\gamma_3+\bar\gamma_3\gamma_4) \\
&\quad +(1-|\gamma_1|^2)(1-|\gamma_2|^2)\gamma_3
(3\bar\gamma_1^2-\bar\gamma_1\gamma_3+3\bar\gamma_1|\gamma_2|^2\gamma_3+\bar\gamma_2^2\gamma_3^2) \\
&\quad -(1-|\gamma_1|^2)\bar\gamma_1^3\gamma_2^4.
\end{align*}

As the reader easily guesses, the following, known as Schur's recurrence relation
({\it cf.}~\cite{Simon:OP1}), is verified by a simple induction argument.

\begin{lem}\label{lem:F}
For each $n\ge 2,$ there exists a function $G_n(\gamma_1,\dots,\gamma_{n-1})$
of $n-1$ complex variables $\gamma_1,\dots,\gamma_{n-1}$ such that the following
equality holds:
$$
F_n(\gamma_1,\dots,\gamma_n)=(1-|\gamma_1|^2)\cdots(1-|\gamma_{n-1}|^2)\gamma_n
+G_n(\gamma_1,\dots,\gamma_{n-1}).
$$
\end{lem}

The following result will be the basis of our arguments below.
It is not new (see, for instance, (1.3.47) in \cite{Simon:OP1}) but a proof is given
for convenience of the reader.

\begin{lem}\label{lem:rec}
Let $\omega(z)=c_1z+c_2z^2+\cdots$ be a function in $\Schur_0$
with its Schur vector $(0, \gamma_1, \gamma_2, \dots).$
Then
$c_n=F_n(\gamma_1,\gamma_2,\dots, \gamma_n)$
for $n=1, 2, 3, \dots.$
\end{lem}

\begin{pf}
We show the assertion by induction on $n.$
When $n=1,$ the assertion is trivial because $c_1=\gamma_1.$
Suppose that the assertion is valid up to $n-1$ for some $n\ge 2.$
Then
\begin{equation}\label{eq:ck}
c_k=F_k(\gamma_1,\gamma_2,\dots, \gamma_k),  \qquad k=1,2,n-1.
\end{equation}
Let $f(z)=\omega(z)/z.$ Then $f\in\Schur$ and it has the Schur vector
$(\gamma_1,\gamma_2,\dots).$
Put $\tilde \omega(z)=z(\sigma f)(z)=zf_1(z)=c_1'z+c_2'z^2+\cdots,$ whose Schur vector is
$(0, \gamma_2, \gamma_3,\dots).$
Therefore, by assumption of the induction again,
\begin{equation}\label{eq:ck'}
c_k'=F_k(\gamma_2,\gamma_3,\dots, \gamma_{k+1}), \qquad k=1,2,n-1.
\end{equation}
By construction, we have the relation
$$
\frac{f(z)-\gamma_1}z=f_1(z)\left(1-\bar\gamma_1 f(z)\right).
$$
We now substitute the Taylor expansions of $f(z)$ and $f_1(z)$
and compare the coefficient of $z^{n-2}$ for both sides to obtain the relation
$$
c_{n}=
(1-|\gamma_1|^2)c_{n-1}'-\bar\gamma_1\sum_{k=2}^{n-1} c_{n-k}'c_k.
$$
Using \eqref{eq:ck} and \eqref{eq:ck'}, we obtain the required formula so that
the induction argument has been completed.
\end{pf}

We remark that for a general function $f(z)=c_0+c_1z+c_2z^2+\cdots$ in $\Schur$ with
the Schur vector $(\gamma_0,\gamma_1,\dots),$ the relations
$$
c_n=F_{n+1}(\gamma_0,\gamma_1,\dots,\gamma_n), \qquad n=0,1,2,\dots,
$$
follow from Lemma \ref{lem:rec}
because the function $zf(z)=c_0z+c_1z^2+\cdots$ belongs to $\Schur_0$
and has $(0,\gamma_0,\gamma_1,\dots)$ as its Schur vector.

\section{Relationship between Schur and Carath\'eodory classes}

As is well known, a function $g\in\PP$ corresponds, in a one-to-one manner,
to a function $\omega\in\Schur_0$ through the Cayley transformation:
\begin{equation}\label{eq:Cayley}
g(z)=\frac{1+\omega(z)}{1-\omega(z)},\quad z\in\D.
\end{equation}
Since it is more natural to think about the quantity $p_n/2$ for
Carath\'eodory functions, we expand $g(z)$ and $\omega(z)$ in the forms
$$
g(z)=1+2\sum_{n=1}^\infty b_nz^n
\aand
\omega(z)=\sum_{n=1}^\infty c_nz^n.
$$
In order to describe $b_n$'s in terms of $c_k$'s, we 
define a sequence of polynomials $Q_n(x_1,\dots, x_n)~(n=1,2,3,\dots)$ 
in the variables $x_1,\dots, x_n$ inductively by $Q_1(x_1)=x_1$ and
\begin{equation}\label{eq:Q}
Q_n(x_1,\dots, x_n)=x_n+\sum_{k=1}^{n-1}x_{n-k} Q_k(x_1,\dots, x_k),
\quad n\ge 2.
\end{equation}
It is worth noting that $Q_n(x_1, \dots, x_n)$ is a polynomial in $x_1,\dots, x_n$
with non-negative integer coefficients.
Then we have the following result.

\begin{lem}\label{lem:bn}
$b_n=Q_n(c_1,\dots, c_n)$ for $n=1,2,\dots.$
\end{lem}

\begin{pf}
By \eqref{eq:Cayley}, we have the relation $(g(z)+1)\omega(z)=g(z)-1.$
We now substitute the above Taylor expansions of $g(z)$ and $\omega(z)$ to obtain
\begin{equation}\label{eq:bc}
b_{n}=\sum_{k=0}^{n-1} c_{n-k}b_k
=c_n+\sum_{k=1}^{n-1} c_{n-k}b_k,
\end{equation}
where we set $b_0=1$ for convenience.
We show the assertion by induction on $n.$
When $n=1$ the assertion is clear.
Suppose the assertion is valid up to $n-1;$ that is,
$b_k=Q_k(c_1,\dots, c_k)$ for $k=1,\dots, n-1.$
Then by the above relation and \eqref{eq:Q}
$$
b_n=c_n+\sum_{k=1}^{n-1} c_{n-k}Q_k(c_1,\dots, c_k)
=Q_n(c_1,\dots, c_n),
$$
which means that the assertion is valid for $n.$
Thus the proof is complete by the induction.
\end{pf}

For instance we have the formulae
\begin{align*}
b_1&=c_1, \\
b_2&=c_1^2+c_2, \\
b_3&=c_1^3+2c_1c_2+c_3, \\
b_4&=c_1^4+3c_1^2c_2+c_2^2+2c_1c_3+c_4, \\
b_5&=c_1^5+4c_1^3c_2+3c_1^2c_3+2c_2c_3+3c_1c_2^2+2c_1c_4+c_5.
\end{align*}

By \eqref{eq:bc}, we have also
$$
c_n=b_n-\sum_{k=1}^{n-1} c_{n-k}b_k,
=b_n-\sum_{k=1}^{n-1} b_{n-k}c_k.
$$
Thus, if we define a sequence of polynomials $R_n(x_1,\dots, x_n)~(n=1,2,3,\dots)$ 
in the variables $x_1,\dots, x_n$ inductively by $R_1(x_1)=x_1$ and
\begin{equation}\label{R}
R_n(x_1,\dots, x_n)=x_n-\sum_{k=1}^{n-1}x_{n-k} R_k(x_1,\dots, x_k),
\qquad n\ge 2,
\end{equation}
then we obtain the following formula in the same way.

\begin{lem}\label{lem:cn}
$c_n=R_n(b_1,\dots, b_n)$ for $n=1,2,\dots.$
\end{lem}

We note that the polynomial mapping $\vec{Q}:\C^n\to\C^n$
It might be interesting to observe that the polynomials $Q_n$ and $R_n$ are related
by a very simple relation.

\begin{prop}
$R_n(x_1,\dots, x_n)=-Q_n(-x_1,\dots, -x_n)$ for $n\ge 1.$
\end{prop}

\begin{pf}
We again use the induction on $n.$
When $n=1,$ the assertion is clear.
Suppose that the assertion is valid up to $n-1.$
Then by definition and the induction assumption we see that
$$
R_n(x_1,\dots, x_n)
=x_n+\sum_{k=1}^{n-1}x_{n-k} Q_k(-x_1,\dots, -x_k)
=-Q_n(-x_1,\dots, -x_n),
$$
which shows the assertion for $n.$
\end{pf}

Therefore, we obtain the following formulae easily from the previous ones:
\begin{align*}
c_1&=b_1, \\
c_2&=-b_1^2+b_2, \\
c_3&=b_1^3-2b_1b_2+b_3, \\
c_4&=-b_1^4+3b_1^2b_2-b_2^2-2b_1b_3+b_4, \\
c_5&=b_1^5-4b_1^3b_2+3b_1^2b_3-2b_2b_3+3b_1b_2^2-2b_1b_4+b_5.
\end{align*}

We end the section with a simple observation.
We define mappings $\vec{Q}_n:\C^n\to\C^n$ and $\vec{R}_n:\C^n\to\C^n$ by
\begin{align*}
\vec{Q}_n(x_1,\dots, x_n)&=(Q_1(x_1), Q_2(x_1,x_2),\dots, Q_n(x_1,\dots, x_n)), \\
\vec{R}_n(x_1,\dots, x_n)&=(R_1(x_1), R_2(x_1,x_2),\dots, R_n(x_1,\dots, x_n)). \\
\end{align*}
Then we get the following result by construction.

\begin{lem}\label{lem:QR}
The mappings $\vec{Q}_n$ and $\vec{R}_n$ are both polynomial automorphisms
of $\C^n$ and they are inverses to each other; namely,
$\vec{Q}_n\circ\vec{R}_n=\vec{R}_n\circ\vec{Q}_n=\id_{\C^n}.$
\end{lem}

\section{Main results}

We now define a sequence of functions $T_n=T_n(\gamma_1,\dots, \gamma_n)$ of
$n$ complex variables $\gamma_1,\dots,\gamma_n$ by
\begin{equation}\label{eq:Tn}
T_n(\gamma_1,\dots, \gamma_n)
=Q_n(F_1(\gamma_1), F_2(\gamma_1,\gamma_2),\dots, F_{n}(\gamma_1,\dots,\gamma_n)),
\end{equation}
where $F_k(\gamma_1,\dots,\gamma_k)$ and $Q_n(x_1,\dots, x_n)$ are defined by
\eqref{eq:F} and \eqref{eq:Q} respectively.
For instance,
\begin{align*}
T_1&=\gamma_1, \\
T_2&=\gamma_1^2+\gamma_2(1-|\gamma_1|^2), \\
T_3&=\gamma_1^3+(1-|\gamma_1|^2)\gamma_2(2\gamma_1-\bar\gamma_1\gamma_2)
+(1-|\gamma_1|^2)(1-|\gamma_2|^2)\gamma_3, \\
T_4&=\gamma_1^4
+(1-|\gamma_1|^2)\gamma_2(3\gamma_1^2-2\gamma_2+\bar\gamma_1^2\gamma_2^2)
+3(1-|\gamma_1|^2)^2\gamma_2^2 \\
&\quad
+(1-|\gamma_1|^2)(1-|\gamma_2|^2)\gamma_3(2\gamma_1-2\bar\gamma_1\gamma_2
-\bar\gamma_2\gamma_3) \\
&\quad +(1-|\gamma_1|^2)(1-|\gamma_2|^2)(1-|\gamma_3|^2)\gamma_4, \\
T_5&=\gamma_1^5
+(1-|\gamma_1|^2)\gamma_2
(4\gamma_1^3-3\gamma_1\gamma_2
+2\bar\gamma_1\gamma_2^2-\bar\gamma_1^3\gamma_2^3) \\
&\quad
+2(1-|\gamma_1|^2)^2\gamma_2^2(3\gamma_1-2\bar\gamma_1\gamma_2) \\
&\quad
+6(1-|\gamma_1|^2)^2(1-|\gamma_2|^2)\gamma_2\gamma_3
-3(1-|\gamma_1|^2)(1-|\gamma_2|^2)^2\bar\gamma_1\gamma_3^2 \\
&\quad
+(1-|\gamma_1|^2)(1-|\gamma_2|^2)\gamma_3
(3\gamma_1^2+3\bar\gamma_1^2\gamma_2^2+2\bar\gamma_1\gamma_3
-2\gamma_1\bar\gamma_2\gamma_3-4\gamma_2+\bar\gamma_2^2\gamma_3^2) \\
&\quad +(1-|\gamma_1|^2)(1-|\gamma_2|^2)(1-|\gamma_3|^2)\gamma_4(2\gamma_1
-2\bar\gamma_1\gamma_2-2\bar\gamma_2\gamma_3-\bar\gamma_3\gamma_4) \\
&\quad +(1-|\gamma_1|^2)(1-|\gamma_2|^2)(1-|\gamma_3|^2)(1-|\gamma_4|^2)\gamma_5.
\end{align*}

Note that the formulae for $T_1, T_2$ and $T_3$ appear as (1.3.51-53) in \cite{Simon:OP1}.
To formulate our main result on coefficients,
it is convenient to consider the coefficient body of order $n$ for a subclass $\F$ of $\hol(\D):$ 
$$
\X_n(\F)=\{(a_0,a_1,\dots, a_n): f(z)=a_0+a_1z+\cdots+a_nz^n+O(z^{n+1})
~\text{for some}~ f\in\F\}.
$$
We remark that $\X_n(\F)$ is convex whenever $\F$ is a convex subset of $\hol(\D).$
%In particular, $\X_n(\PP)$ is convex for every $n\ge 1.$
%Trivial cases are $\X_0(\Schur)=\bD$ and $\X_1(\PP)=\{(1,a):|a|\le 2\}.$
We recall that a subset $A$ of $\R^d$ is called a {\it convex body} if
$A$ is compact and convex and has non-empty interior.
It is well known that a convex body in $\R^d$ is homeomorphic to the
closed unit ball $\mathbb{B}^{d}=\{(x_1,\cdots, x_{d})\in\R^{d}: x_1^2+\cdots+x_{d}^2\le 1\}$ (see \cite[\S 11.3]{Berger:geom1}).

\begin{thm}\label{thm:main}
Let $n$ be a positive integer.
The coefficient body $\X_n(\PP)$ of order $n$ for the Carath\'eodory class $\PP$ is expressed
as $\{1\}\times 2\V_n,$ where $\V_n$ is a convex body
in $\C^n.$
Moreover, 
$$
\vec{T}_n(\gamma_1,\dots, \gamma_n)=(T_1(\gamma_1), T_2(\gamma_1,\gamma_2),\dots,
T_n(\gamma_1,\dots,\gamma_n)).
$$
is a continuous mapping of $\bD^n$ onto $\V_n$
and satisfies $\vec{T}_n(\D^n)=\Int\V_n$ and $\vec{T}_n(\partial\bD^n)=\partial\V_n.$
In addition,
$\vec{T}_n:\D^n\to\Int\V_n$ is a real analytic diffeomorphism 
but $\vec{T}_n$ is not injective on the boundary $\partial\bD^n$ of $\bD^n$ for $n=2,3,\dots.$
\end{thm}

\begin{pf}
By the normalization condition for Carath\'eodory functions, we first observe that
$\X_n(\PP)$ can be expressed as the form $\{1\}\times 2\V_n,$ where
$V_n=\{(b_1,\dots,b_n)\in\C^n: (1,2b_1,\dots, 2b_n)\in \X_n(\PP)\}.$
Since $\PP$ is a convex subset of $\hol(\D),$ it is evident that $\V_n$ is convex
in $\C^n.$
Similarly, we can write $\X_n(\Schur_0)=\{0\}\times \U_n.$

As we saw in the previous section, the coefficients of a function $g(z)=1+2b_1z+2b_2z^2+\cdots$
in $\PP$ and those of $\omega(z)=c_1z+c_2z^2+\cdots$ in $\Schur_0$ are related by
$(b_1,\dots, b_n)=\vec{Q}_n(c_1,\dots, c_n)$ whenever $g$ and $\omega$ are
related by $g=(1+\omega)/(1-\omega).$
In particular, we have the relation $\V_n=\vec{Q}_n(\U_n).$
Let
$$
\vec{F}_n(\gamma_1,\dots, \gamma_n)=(F_1(\gamma_1), F_2(\gamma_1,\gamma_2),\dots,
F_n(\gamma_1,\dots,\gamma_n)),
$$
where $F_1,\dots, F_n$ are defined in \eqref{eq:F}.
Then $\vec{T}_n=\vec{Q}_n\circ\vec{F}_n$ by construction.
Since $\vec{Q}_n:\C^n\to\C^n$ is a polynomial automorphism of $\C^n$ by Lemma \ref{lem:QR}, 
the other assertions follow from the next proposition,
which may be of independent interest.
\end{pf}

\begin{prop}\label{prop:diffeo}
Let $n$ be a positive integer.
The coefficient body $\X_n(\Schur_0)$ of order $n$ for $\Schur_0$ is described by
$\X_n(\Schur_0)=\{0\}\times \U_n,$ where $\U_n$ is a convex body in $\C^n.$
Moreover, $\vec{F}_n$ maps $\bD^n$ continuously onto $\U_n$
and satisfies $\vec{F}_n(\D^n)=\Int\U_n$ and $\vec{F}_n(\partial\bD^n)=\partial\U_n.$
Furthermore, $\vec{F}_n:\D^n\to\Int\U_n$ is a real analytic diffeomorphism 
but $\vec{F}_n$ is not injective on $\partial\bD^n$ for $n=2,3,\dots.$
\end{prop}

Before the proof, we make a preliminary observation.

\begin{lem}\label{lem:inj}
Let $(\gamma_1,\dots,\gamma_n)\in\D^{n-1}\times\C$ and $(\gamma_1',\dots,\gamma_n')\in\C^n.$
If $\vec{F}_n(\gamma_1,\dots,\gamma_n)=\vec{F}_n(\gamma_1',\dots,\gamma_n'),$
then $(\gamma_1,\dots,\gamma_n)=(\gamma_1',\dots,\gamma_n').$
\end{lem}

\begin{pf}
We show $\gamma_j=\gamma_j'$ by induction.
Since $\gamma_1=F_1(\gamma_1)=F_1(\gamma_1')=\gamma_1',$
the assertion is clear for $j=1.$
Assume next that the assertion holds true up to $j-1;$
that is $\gamma_1=\gamma_1',\dots, \gamma_{j-1}=\gamma_{j-1}'.$
Then, by Lemma \ref{lem:F}, we observe
\begin{align*}
& \quad \
(1-|\gamma_1|^2)\cdots(1-|\gamma_{j-1}|^2)\gamma_j+G_j(\gamma_1,\dots,\gamma_{j-1}) \\
&=F_j(\gamma_1,\dots,\gamma_j) \\
&=F_j(\gamma_1',\dots,\gamma_j') \\
&=(1-|\gamma_1'|^2)\cdots(1-|\gamma_{j-1}'|^2)\gamma_j'+G_j(\gamma_1',\dots,\gamma_{j-1}') \\
&=(1-|\gamma_1|^2)\cdots(1-|\gamma_{j-1}|^2)\gamma_j'+G_j(\gamma_1,\dots,\gamma_{j-1}).
\end{align*}
Since $(1-|\gamma_1|^2)\cdots(1-|\gamma_{j-1}|^2)\ne0$ by assumption,
we conclude that $\gamma_j=\gamma_j'.$
Thus the induction argument has been completed.
\end{pf}

We are ready to show the above proposition.

\begin{pf}[Proof of Proposition \ref{prop:diffeo}]
Since $\Schur_0$ is compact and convex in the vector space $\hol(\D)$ endowed with the topology
of locally uniform convergence on compact subsets of $\D,$
$\X_n(\PP)$ is also compact and convex.
Consequently, the set $\U_n$ is compact and convex in $\C^n.$
In order to see that $\U_n$ is a convex body,
we have only to show that $\U_n$ has a non-empty interior.
For instance, we see that the polynomial $\omega(z)=c_1z+c_2z^2+\cdots+c_nz^n$
with $|c_1|+|c_2|+\cdots +|c_n|<1$ is contained in the class $\Schur_0.$
Hence, the non-empty open set $|c_1|+|c_2|+\cdots+|c_n|<1$ is contained in $\V_n.$

We next show that $\vec{F}_n(\bD^n)=\U_n.$
Indeed, for $(c_1,\dots, c_n)\in\U_n,$ by definition, there is a function $\omega\in\Schur_0$
such that $\omega(z)=c_1z+\cdots+c_nz^n+O(z^{n+1}).$
Let $(0,\gamma_1,\gamma_2,\dots)$ be the Schur vector of $\omega.$
Then $(c_1,\dots,c_n)=\vec{F}_n(\gamma_1,\dots,\gamma_n).$
Thus the inclusion relation $\vec{F}_n(\bD^n)\supset\U_n$ follows.
Conversely, choose $n$ points $\gamma_1,\dots,\gamma_n$ from $\bD$ arbitrarily.
When $(\gamma_1,\dots,\gamma_n)\in\D^n,$ let $\vec{\gamma}=(\gamma_1,\dots,\gamma_n,0,0,\dots).$
Otherwise, there is a unique $m$ such that $|\gamma_1|<1,\dots, |\gamma_{m-1}|<1$ and
$|\gamma_m|=1.$
Then we set $\vec{\gamma}=(\gamma_1,\dots,\gamma_m,0,0,\dots).$
Theorem C now guarantees existence of a function $f$ such that
$\vec{\gamma}(f)=\vec{\gamma}.$
We can now expand the function $\omega(z)=zf(z)$ 
in the form $\omega(z)=c_1z+c_2z^2+\cdots.$
Thus we see that $\vec{F}_n(\gamma_1,\dots,\gamma_n)=(c_1,\dots,c_n)\in\U_n.$
Hence, we have shown that the other inclusion relation $\vec{F}_n(\bD^n)\subset\U_n.$
We have shown $\vec{F}_n(\bD^n)=\U_n.$

Note that $\vec{F}_n$ is real analytic on $\C^n.$
We show now that $\vec{F}_n:\D^n\to\C^n$ is locally diffeomorphic; in other words,
the Jacobian $J_{\vec{F}_n}$ does not vanish on $\D^n.$
For a moment, we set $\delta_j=(1-|\gamma_1|^2)\cdots(1-|\gamma_j|^2)$ for short.
Regarding $\vec{F}_n$ as a column vector, with the help of Lemma \ref{lem:F}, we compute
\begin{align*}
J_{\vec{F}_n}&=
\begin{vmatrix}
\dfrac{\partial F_1}{\partial\gamma_1} & \dfrac{\partial F_1}{\partial\bar\gamma_1}& \cdots&
\dfrac{\partial F_1}{\partial\gamma_n} & \dfrac{\partial F_1}{\partial\bar\gamma_n} \medskip \\
\dfrac{\partial \overline{F_1}}{\partial\gamma_1} & 
\dfrac{\partial \overline{F_1}}{\partial\bar\gamma_1}& \cdots&
\dfrac{\partial \overline{F_1}}{\partial\gamma_n} & \dfrac{\partial\overline{F_1}}{\partial\bar\gamma_n} \\
\null & \null & \ddots & \null & \null  \\
\dfrac{\partial F_n}{\partial\gamma_1} & \dfrac{\partial F_n}{\partial\bar\gamma_1}& \cdots&
\dfrac{\partial F_n}{\partial\gamma_n} & \dfrac{\partial F_n}{\partial\bar\gamma_n} \medskip \\
\dfrac{\partial \overline{F_n}}{\partial\gamma_1} & 
\dfrac{\partial \overline{F_n}}{\partial\bar\gamma_1}& \cdots&
\dfrac{\partial \overline{F_n}}{\partial\gamma_n} & \dfrac{\partial\overline{F_n}}{\partial\bar\gamma_n}
\end{vmatrix} 
=\begin{vmatrix}
1 & 0 & 0 & 0 & \cdots & 0 & 0 \\
0 & 1 & 0 & 0 & \cdots & 0 & 0 \\
* & * & \delta_1 & 0 & \cdots & 0 & 0 \\
* & * & 0 & \delta_1 & \cdots & 0 & 0 \\
\null & \null & \null & \null & \ddots & \null & \null  \\
* & * & * & * & \cdots & \delta_{n-1} & 0 \\
* & * & * & * & \cdots & 0 & \delta_{n-1}
\end{vmatrix} \\
&=\delta_1^2\delta_2^2\cdots\delta_{n-1}^2>0.
\end{align*}
Thus we have shown that $\vec{F}_n$ is locally diffeomorphic on the domain $\D^n.$
In particular, the image $\vec{F}_n(\D^n)$ is contained in the interior
$\Int \U_n$ of $\U_n.$
Moreover, Lemma \ref{lem:inj} now implies that $\vec{F}_n$ is injective on $\D^n.$

We next prove that $\vec{F}_n$ maps $\partial\bD^n$ onto $\partial\U_n.$
%Recall that $\vec{F}_n:\C^n\to\C^n$ is real analytic.
Let $(\gamma_1,\dots,\gamma_n)\in\partial\bD^n.$
Then $|\gamma_1|<1,\dots, |\gamma_{j-1}|<1$ and $|\gamma_j|=1$ for some
$1\le j\le n.$
To the contrary, we suppose that $p=\vec{F}_n(\gamma_1,\dots,\gamma_n)\in\Int\U_n.$
First consider the case when $j=n.$
Then for small enough $\delta>0,$ we have
$p_\delta=\vec{F}_n(\gamma_1,\dots,\gamma_{n-1},(1+\delta)\gamma_n)\in\Int\U_n.$
Since $\vec{F}_n:\bD^n\to\U_n$ is surjective, $p_\delta=\vec{F}_n(\gamma_1',\dots,\gamma_n')$
for some $(\gamma_1',\dots,\gamma_n')\in\bD^n.$
We apply Lemma \ref{lem:inj} again to deduce that
$(\gamma_1,\dots,\gamma_{n-1},(1+\delta)\gamma_n)=(\gamma_1',\dots,\gamma_n').$
In particular, $(1+\delta)\gamma_n=\gamma_n'\in\bD,$ which contradicts 
$(1+\delta)|\gamma_n|=1+\delta>1.$
Therefore, this case does not occur.
When $1\le j<n,$ we consider the projection $\pi:\C^n\to\C^j$ defined by
$\pi(z_1,\dots,z_n)=(z_1,\dots, z_j).$
By definition, we have $\pi(\U_n)=\U_j.$
Since $\pi$ is an open mapping, we have $\pi(p)\in\Int \U_j.$
However, since $\pi(p)=\vec{F}_j(\gamma_1,\dots,\gamma_j),$
this is again impossible by the same reason.
Therefore $\vec{F}_n(\gamma_1,\dots,\gamma_n)\in\partial\U_n$ for
$(\gamma_1,\dots,\gamma_n)\in\partial\bD^n$ at any event.
Recalling that $\vec{F}_n:\bD^n\to\U_n$ is surjective,
we now conclude that $\vec{F}_n(\D^n)=\Int\U_n$
and that $\vec{F}_n(\partial\bD^n)=\partial\U_n.$

Finally, we see that $\vec{F}_n$ is not injective on $\partial\bD^n$ when $n\ge 2.$
Indeed, by Lemma \ref{lem:F}, $\vec{F}_n(\gamma_1,\dots,\gamma_{n-1},\gamma_n)
=\vec{F}_n(\gamma_1,\dots,\gamma_{n-1}, 0)$ for any
$(\gamma_1,\dots,\gamma_{n-1})\in\partial\bD^{n-1}$ and $\gamma_n\in\bD.$
\end{pf}

We now make a comparison with the Libera-Z\l otkiewicz lemma (Theorem B above).
In terms of $b_n=p_n/2,$ we can reformulate it as
\begin{align*}
b_2&=b_1^2+x(1-b_1^2)\quad\text{and} \\
b_3&=b_1^3+(1-b_1^2)x(2b_1-b_1x)+(1-b_1^2)(1-|x|^2)y.
\end{align*}
We observe that their results agree with our formulae for $T_2$ and $T_3$
when $b_1=\gamma_1\ge0, x=\gamma_2, y=\gamma_3.$

\section{Recursion for $T_n(\gamma_1,\dots,\gamma_n)$}

In the previous section, we defined $T_n(\gamma_1,\dots,\gamma_n)$ as the
composition of the polynomial $Q_n(x_1,\dots,x_n)$ with $\vec{F}_n(\gamma_1,\dots,\gamma_n).$
Recall that $Q_n$ and $F_k$ are both defined recursively.
In principle, there should be a recursive formula which defines $T_n.$
We end this note by giving such a formula.
Let $g(z)=1+2(b_1z+b_2z^2+\cdots)$ be a function in $\PP$ and 
take $f\in\Schur$ with $\vec{\gamma}(f)=(\gamma_1,\gamma_2,\dots)$
so that $g(z)=(1+zf(z))/(1-zf(z)).$
Then, by construction of $T_n,$ we obtain
\begin{equation}\label{eq:bk}
b_k=T_k(\gamma_1,\dots,\gamma_k), \qquad k=1,2,3,\cdots.
\end{equation}
Let $g_1(z)=(1+zf_1(z))/(1-zf_1(z))=1+2(b_1'z+b_2'z^2+\cdots)$ for 
$f_1(z)=\sigma f(z)=(f(z)-\gamma_1)/[z(1-\bar\gamma_1f(z))].$
Since $\vec{\gamma}(f_1)=(\gamma_2,\gamma_3,\dots),$ we have
\begin{equation}\label{eq:bk'}
b_k'=T_k(\gamma_2,\dots,\gamma_{k+1}), \qquad k=1,2,3,\cdots.
\end{equation}
In view of the relation $zf(z)=(g(z)-1)/(g(z)+1),$ we obtain
$$
\frac{g_1(z)-1}{g_1(z)+1}=zf_1(z)=\frac{f(z)-\gamma_1}{1-\bar\gamma_1f(z)}
=\frac{g(z)-1-\gamma_1z(g(z)+1)}{z(g(z)+1)-\bar\gamma_1(g(z)-1)},
$$
from which we derive the formula
$$
\big\{g_1(z)-1\big\}\big\{z(g(z)+1)-\bar\gamma_1(g(z)-1)\big\}
=\big\{g_1(z)+1\big\}\big\{g(z)-1-\gamma_1z(g(z)+1)\big\}.
$$
Substituting the power series expansions of $g(z)$ and $g_1(z),$ we get the relation
$$
\sum_{n=1}^\infty b_n'z^n\left[\sum_{n=1}^\infty b_{n-1}z^{n}
-\bar\gamma_1\sum_{n=1}^\infty b_nz^n\right]
=\sum_{n=0}^\infty b_n'z^n\left[\sum_{n=1}^\infty b_nz^{n}
-\gamma_1\sum_{n=1}^\infty b_{n-1}z^{n}\right],
$$
where we set $b_0=b_0'=1.$
We look at the coefficients of $z^n$ of the functions in the both sides to obtain
\begin{align*}
\sum_{k=1}^{n-1} b_k'(b_{n-k-1}-\bar\gamma_1b_{n-k})
&=\sum_{k=0}^{n-1} b_k'(b_{n-k}-\gamma_1b_{n-k-1}) \\
&=b_n-\gamma_1b_{n-1}+\sum_{k=1}^{n-1} b_k'(b_{n-k}-\gamma_1b_{n-k-1}).
\end{align*}
Hence,
$$
b_n=\gamma_1b_{n-1}+\sum_{k=1}^{n-1} b_k'\Big[(1+\gamma_1)b_{n-k-1}
-(1+\bar\gamma_1)b_{n-k}\Big].
$$
We now substitute \eqref{eq:bk} and \eqref{eq:bk'} into the last formula to have
the following result.

\begin{thm}
The functions $T_n(\gamma_1,\dots,\gamma_n)$ defined in \eqref{eq:Tn} are
described by the following recursive formula with the initial condition $T_0=1:$
\begin{align*}
&\qquad\qquad T_n(\gamma_1,\dots,\gamma_n)
=\gamma_1 T_{n-1}(\gamma_1,\dots,\gamma_{n-1}) \\
&+\sum_{k=1}^{n-1} T_k(\gamma_2,\dots, \gamma_{k+1})
\Big[(1+\gamma_1)T_{n-k-1}(\gamma_1,\dots, \gamma_{n-k-1})
-(1+\bar\gamma_1)T_{n-k}(\gamma_1,\dots, \gamma_{n-k})\Big].
\end{align*}
\end{thm}

We remark that the transformation $g_1$ from $g$ above was already considered
by Brown \cite{Brown87} in a more general context.

\def\cprime{$'$} \def\cprime{$'$} \def\cprime{$'$}
\providecommand{\bysame}{\leavevmode\hbox to3em{\hrulefill}\thinspace}
\providecommand{\MR}{\relax\ifhmode\unskip\space\fi MR }
\providecommand{\MRhref}[2]{%
  \href{http://www.ams.org/mathscinet-getitem?mr=#1}{#2}
}
\providecommand{\href}[2]{#2}

%\bibliography{papers}
\end{document}